\documentclass{article}

\usepackage{natbib}
\usepackage{mathrsfs}
\usepackage{authblk}
\usepackage{floatrow}
\usepackage{natbib}
\usepackage{mathrsfs}
\usepackage{floatrow}

\usepackage[utf8]{inputenc} 
\usepackage[T1]{fontenc}    
\usepackage{hyperref}       
\usepackage{url}            
\usepackage{booktabs}       
\usepackage{amsfonts}       
\usepackage{nicefrac}       
\usepackage{microtype}      
\usepackage{xcolor}         
\usepackage{microtype}
\usepackage{graphicx}
\usepackage{subfigure}
\usepackage{booktabs} 
\usepackage{amsfonts}

\usepackage{amsmath}
\usepackage{amssymb}
\usepackage{mathtools}
\usepackage{amsthm}

\usepackage[capitalize,noabbrev]{cleveref}

\theoremstyle{plain}

\theoremstyle{definition}

\theoremstyle{remark}

\newtheorem*{theorem-non}{Theorem}
\usepackage[textsize=tiny]{todonotes}

\usepackage{hyperref}

     \usepackage[final]{neurips_2021}

\usepackage[utf8]{inputenc} 
\usepackage[T1]{fontenc}    
\usepackage{hyperref}       
\usepackage{url}            
\usepackage{booktabs}       
\usepackage{amsfonts}       
\usepackage{nicefrac}       
\usepackage{microtype}      
\usepackage{xcolor}         
\usepackage{microtype}
\usepackage{graphicx}
\usepackage{subfigure}
\usepackage{booktabs} 
\usepackage{amsfonts}

\usepackage{amsmath}
\usepackage{amssymb}
\usepackage{mathtools}
\usepackage{amsthm}

\usepackage[capitalize,noabbrev]{cleveref}

\theoremstyle{plain}

\usepackage[textsize=tiny]{todonotes}

\usepackage{hyperref}

\title{Private independence testing across two parties}

\begin{document}

\author[1]{\small Praneeth Vepakomma}

\author[1]{\small Mohammad Mohammadi Amiri}
\author[2]{\small Clément L. Canonne}
\author[1]{\\\small Ramesh Raskar}
\author[1]{\small Alex Pentland}

\affil[1]{\footnotesize Massachusetts Institute of Technology, USA}
\affil[2]{\footnotesize The University of Sydney, Australia}

\maketitle

\begin{abstract}
We introduce $\pi$-\textsf{test}, a privacy-preserving algorithm for testing statistical independence between data distributed across multiple parties. Our algorithm relies on privately estimating the distance correlation between datasets, a quantitative measure of independence introduced in~\cite{szekely2007measuring}. We establish both additive and multiplicative error bounds on the utility of our differentially private test, which we believe will find applications in a variety of distributed hypothesis testing settings involving sensitive data.
\end{abstract}

\section{Introduction}
Hypothesis testing is a staple of theoretical and applied statistics, whose impact in decision making has gone beyond that of statistics into a wide-range of adjacent and seemingly far-flung fields of science and engineering. At a high level, hypothesis testing formalizes the choice to reject (or not reject) a \emph{null} hypothesis $\mathcal{H}_{0}$ in comparison to an alternate hypothesis $\mathcal{H}_{1}$, where both hypotheses are defined as disjoint sets of probability distributions. The decision is to be made with a high probability upon looking at samples from some unknown distribution, assumed to belong to $\mathcal{H}_{0} \cup \mathcal{H}_{1}$.

Examples of hypothesis testing abound, including, to name a few, goodness-of-fit and testing validity of a purported statistical model, testing stationarity of a time series, hypothesis tests for sequential decision making such as change-point detection, and so on. Central to our study is one of the most fundamental hypotheses testing problems, dating back to~\citet{Pearson1900}: that of \emph{statistical independence}. In contexts where the underlying data (i.e., samples) is personal or conveys sensitive information, guaranteeing privacy (of data samples) in hypothesis testing has recently become the focus of much theoretical and practical interest, both for legal and societal reasons. Indeed, privacy-preserving computation is necessitated by a range of motivations, including regulatory compliance, personal preferences, confidentiality, as well a a migration to distributed forms of computation where client data or intermediate computations of their data are shared among various (not necessarily trusted) parties. In this paper, we quantify privacy in the framework of \emph{differential privacy}, a rigorous and widely used notion of privacy introduced in~\citet{dwork2014algorithmic}, and focus on a two-party setting.  
The following is the main problem of consideration in this paper.

\textbf{Problem statement:} How can the sample \emph{test-statistic} for a hypothesis test of independence between two random variables of arbitrary dimension be estimated \emph{privately} when samples from these random
variables are held by two different parties? The central importance of a test-statistic (to hypothesis testing) is detailed in the preliminaries (Section~\ref{prelim}).

\begin{figure}[!htbp]
    \centering
    \includegraphics[scale=0.75]{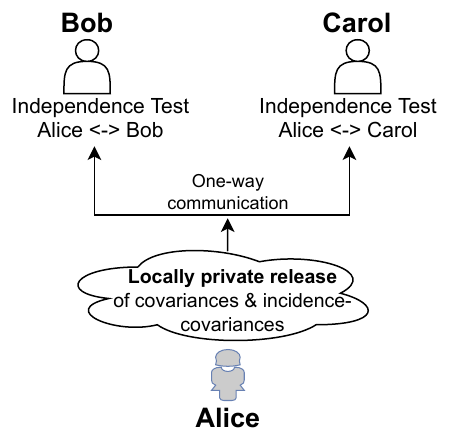}
    \caption{Model of privacy considered, where Alice releases a private summary of data as a sufficient statistic which allows any analyst (here, either Bob or Carol) to test independence with another data set they hold with respect to that of Alice's original data. This private statistic suffices to compute the rest of the test statistic for independence testing on any analyst's premise. We provide utility bounds on the computation of this test-statistic in Section~\ref{utilBds}.}
    \label{fig:setup}
\end{figure}

\textbf{Communication setup:} In terms of communication, we consider a two-party setup, with Alice holding $\boldsymbol{X}$ and Bob holding $\boldsymbol{Y}$. The communication is one-way, from Alice to Bob: Alice sends an intermediate computation in a privatized form to Bob, who then finishes the rest of the computation on its premise. Bob never reveals the answer to Alice and therefore the on-device computation of Bob involves its own non-privatized data. Therefore, Bob obtains the results (on its side) of the hypothesis test of independence between its data and that of Alice's, while maintaining the privacy of Alice's (already privatized) data, regardless of the one-way communications that happen. We refer to this setup as \emph{one-way local privacy} or \emph{one-way local differential privacy}, as the privatization mechanism used by Alice ensures differential privacy \cite{dwork2014algorithmic}. A more general way to think about it is as a data summary that when published by Alice, allows any analyst to test independence with another data set they themselves hold. We illustrate this setup in Figure~\ref{fig:setup}.

\vspace{-3mm}\subsection{Related Work} 
In view of the large amount of literature addressing the problem of independence testing over the years in a variety of statistical settings, we focus here on the closest and most relevant to our work.
\noindent\textbf{Private independence testing for contingency tables.} This first line of work, which aims to develop differentially private independence tests, itself comes in two distinct flavors: one is the so-called \emph{asymptotic regime} (limiting distribution and properties of the test statistic as sample size goes to infinity), and the second is the \emph{finite-sample regime} (coarser utility guarantees, but with explicit, finite bounds on the sample size required to achieve them).
\begin{enumerate}
 \item{\it{}Asymptotic tests:} The work of~\citet{gaboardi2016differentially} considers independence testing in the (central) model of differentially privacy, from an \emph{asymptotic} perspective: namely, they propose a differentially private analogue of the classical chi-squared tests of independence, and analyse the limiting distribution of the test-statistic along with its resulting power (1-Type II error). Note that, Type II error is the probability of failing to reject null hypothesis when the null hypothesis is not true.
              \item{\it{}Non-asymptotic tests:} Later work by~\citet{sheffet2018locally} focuses on the \emph{finite sample} (non-asymptotic) version of the test under a more stringent model of local privacy. The question is formulated, in a manner that is standard in theoretical computer science (specifically in the domain of distribution testing) and in minimax analysis in statistics as part of a composite hypthesis testing problem. They consider a set of product distributions as part of the null hypothesis and frame the alternative hypothesis to contain all the distributions that are ``far'' from product distributions, where ``farness'' is quantified by the total variation distance. They focus on the minimax sample complexity achievable using a specific locally private mechanism of Randomized Response~\citep{warner1965randomized}. Subsequent work by \citet{acharya2021inference} improves on these results, by showing a way to achieve significantly lower sample complexity (still in the locally private setting) by considering a different privacy mechanism and upon establishing matching lower bounds.
\end{enumerate}
    
\noindent\textbf{Non-private independence testing with dependency measures.}     
There have been advances in development of various statistical dependency measures, and an active route of modern independence testing is based on derivation of test-statistics that depend on these measures in addition to on other required terms. We now share some related works that fall in this category. Distance covariance was introduced in \citep{szekely2007measuring} and can be expressed as a weighted $L_{2}$ norm between the characteristic function of the joint distribution and the product of the marginal characteristic functions. This concept has also been studied in high dimensions \citep{szekely2013distance,yao2018testing}, and for testing independence of several random vectors \citep{gao2021asymptotic}. In \citep{sejdinovic2013equivalence}, tests based on distance covariance were shown to be equivalent to a reproducing kernel Hilbert space (RKHS) test for a specific choice of kernel. RKHS tests have been widely studied in the machine learning community, with a survey of the subject given by \citep{harchaoui2013kernel} and \citep{gretton2005measuring,heller2016consistent} in which the Hilbert-Schmidt independence criterion was proposed. These tests are based on embedding the joint distribution and product of the marginal distributions into a Hilbert space and considering the norm of their difference in this space. One drawback of the kernel paradigm here is the computational complexity, though \citep{jitkrittum2017linear,jitkrittum2016adaptive} and \citep{zhang2018large} have recently attempted to address this issue. 

\noindent\textbf{Conditional independence tests and causal discovery.} A conditional measure of dependence called conditional distance correlation was introduced by \citep{wang2015conditional}. The works in \citep{zhang2012kernel,zhang2017causal,wang2020towards,shen2022chi} performed conditional independence testing and applied them to the problem of causal discovery.

 To the best of our knowledge, no other work addresses the question of independence testing under differential privacy (be it local or central) other than the differentially private distribution estimation approach (under total variation distance)~\citep{DiakonikolasHS15} or two-sample goodness-of-fit~\citep{AcharyaSZ18,AliakbarpourDR18} that can be used to obtain sub-optimal sample complexity guarantees for this problem. We note that this body of work differs from ours, both in the distributed model assumed (the way the data is partitioned across users) and in the guarantee provided (dependency measure used). They also restrict themselves to the discrete setting as opposed to ours. Our method can also be applied to test between samples lying in different dimensions. In particular, most of the works discussed above (except for~\citet{gaboardi2016differentially}), follow the norm in distribution testing and focus on very stringent notion of minimax testing under total variation distance, which might be overly conservative in many settings. 
\subsection{Our contributions}
Our contributions are threefold, and can be summarized as follows:
\begin{enumerate}
    \item We provide a mechanism to privatize the test-statistic required to perform a hypothesis test of independence between $\boldsymbol{x} \in \mathbb{R}^d$ hosted by Alice and $\boldsymbol{y} \in \mathbb{R}^m$ hosted by Bob under the regime of one-way local differential privacy.
    \item We show one can express our test statistic as a ratio of sums of \emph{directional variance} queries (of non-private data) with respect to `specific' covariance matrices of sensitive data (for the formal definition of \emph{directional variance} queries, the reader is referred to Section~\ref{prelim}). This reduction enables us to privatize the test-statistic via privatization of these covariances.  This brings up the next question of deriving utility guarantees for the privatized test-statistic, leading to our next contribution. 
    \item We derive both lower and upper bounds on the utility of our private estimator in terms of additive and multiplicative errors. For a chosen $0<\eta<1$, to achieve $(\epsilon,\delta)$-differential privacy our approach results in a multiplicative error factor of $1\pm\eta$, and an additive error $\Delta$ bounded as \[
    -\frac{(1-\eta)^{2}}{2(1+\eta)}\leq \Delta \leq \frac{2 \tau}{(1-\eta)[(1-\eta) s-\tau]},
    \]
    where 
    \[
    \tau \coloneqq \left(\frac{2048 \ln (2 /(m+n) \nu) \ln (2 / \delta)}{\eta \epsilon^{2}}\right) \ln ^{2}\!\left(\frac{128 \ln\frac{1}{(m+n) \nu}}{\eta^{2} \delta}\right),
    \] 
    for some $s>\frac{\tau}{1-\eta}$. The additive and multiplicative approximations hold together at the same time with probability at least $1-(m+n) \nu$ for the user's choice of $\nu$.  
    
\end{enumerate}
\bibliography{example_paper}

\newpage

\end{document}


\maketitle
\appendix
\section{Lower bound} \label{lowerProof}
\begin{theorem}
For $\hat{S}>\frac{n \tau}{1-n}$ and $n \gg m$, we have the following lower bound on $\frac{\bar{\Omega}^2}{\bar{S}}$.
$$
\mathbb{P}\left( \frac{\bar{\Omega}^2}{\bar{S}} \ge \frac{1-\eta}{1+\eta}\left(\frac{\hat{\Omega}}{\hat{S}}\right)- \frac{(1-\eta)^2}{2(1+\eta)}  \right) \geq 1-(m+n) \nu.
$$
\end{theorem}
\begin{proof}From the naive re-arrangement of the individual bounds, we have the following with probability $\geq 1- (m+n) \nu$ 
\begin{align}\label{naive}
\frac{(1-\eta)\hat{\Omega}^2 -m\tau}{(1+\eta)\hat{S} +n \tau} \le \frac{\bar{\Omega}^2}{\bar{S}} \le \frac{(1+\eta)\hat{\Omega}^2 +m \tau}{(1-\eta)\hat{S} -n \tau}.
\end{align}
Rearranging the lower bound on $\frac{\bar{\Omega}^2}{\bar{S}}$ yields
\begin{align}\label{LowerBound1}
\frac{\bar{\Omega}^2}{\bar{S}} & \geqslant \frac{(1-\eta) \hat{\Omega}^2-m \tau}{(1+\eta) \hat{S}+n\tau}= \frac{(1-\eta) \hat{\Omega}^2-m\tau}{(1+\eta) \hat{S}\left(1+ \frac{n \tau}{(1+\eta) \hat{S}}\right)}=\frac{(1-\eta) \hat{\Omega}^2 \left(1+\frac{n\tau}{(1+\eta) \hat{S}}\right)-m \tau-\frac{(1-\eta) \hat{\Omega}^2 n \tau}{(1+\eta) \hat{S}}}{(1+\eta) \hat{S}\left(1+ \frac{n \tau}{(1+\eta) \hat{S}}\right)}\nonumber\\
&=\left(\frac{1-\eta}{1+\eta}\right) \frac{\hat{\Omega}^2}{\hat{S}} - \frac{m \tau (1+\eta) \hat{S} + n \tau (1-\eta) \hat{\Omega}^2}{(1+\eta) \hat{S} \left((1+\eta) \hat{S}+n \tau \right)}.
\end{align}
For $\hat{S} > \frac{n \tau}{1-\eta}$, we have
\newcommand\firstinequal{\mathrel{\overset{\makebox[0pt]{\mbox{\normalfont\tiny\sffamily (a)}}}{\le}}}
\newcommand\secondinequal{\mathrel{\overset{\makebox[0pt]{\mbox{\normalfont\tiny\sffamily (b)}}}{\le}}}
\begin{align}\label{LowerBound2}
\frac{m \tau (1+\eta) \hat{S} + n \tau (1-\eta) \hat{\Omega}^2}{(1+\eta) \hat{S} \left((1+\eta) \hat{S}+n \tau \right)} \firstinequal \frac{m \tau (1+\eta) + n \tau (1-\eta)}{(1+\eta) \left((1+\eta) \hat{S}+n \tau \right)} \secondinequal \frac{m (1+\eta) + n (1-\eta)}{2 n \left(\frac{1+\eta}{1-\eta}\right)},   
\end{align}
where (a) follows from Lemma \ref{boundOmegaS2}, i.e., $\hat{\Omega}^2 \le \hat{S}$, and (b) follows by replacing $\hat{S}$ with $\frac{n \tau}{1-\eta}$ since $\hat{S} > \frac{n \tau}{1-\eta}$.
For $n \gg m$, we can approximate \eqref{LowerBound2} as follows 
\begin{align}\label{LowerBoundLowerApprox}
\frac{m \tau (1+\eta) \hat{S} + n \tau (1-\eta) \hat{\Omega}^2}{(1+\eta) \hat{S} \left((1+\eta) \hat{S}+n \tau \right)} \le \frac{m (1+\eta) + n (1-\eta)}{2 n \left(\frac{1+\eta}{1-\eta}\right)} \approx \frac{(1-\eta)^2}{2(1+\eta)}.   
\end{align}
Plugging \eqref{LowerBoundLowerApprox} into \eqref{LowerBound1} yields 
\begin{align}
\frac{\bar{\Omega}^2}{\bar{S}} \ge \left(\frac{1-\eta}{1+\eta}\right) \frac{\hat{\Omega}^2}{\hat{S}} - \frac{(1-\eta)^2}{2(1+\eta)}   
\end{align}


\end{proof}

\section{Upper bound}
\begin{theorem}
For some $s>\frac{\tau}{1-\eta}$ and $n \gg m$, we have the following upper bound on $\frac{\bar{\Omega}^2}{\bar{S}}$.
$$
\mathbb{P}\left( \frac{\bar{\Omega}^2}{\bar{S}} \leqslant\frac{1+\eta}{1-\eta}\left(\frac{\hat{\Omega}^2}{\hat{S}}\right) + \frac{\tau}{(1-\eta) s-\tau}\right) \geq 1-(m+n) \nu.
$$
\end{theorem}
\label{upperProof}
\begin{proof}
Rearranging the upper bound on $\frac{\bar{\Omega}^2}{\bar{S}}$, given in \eqref{naive}, yields
\begin{align}\label{UpperBound1}
\frac{\bar{\Omega}^2}{\bar{S}} & \le \frac{(1+\eta) \hat{\Omega}^2+m\tau}{(1-\eta) \hat{S}\left(1- \frac{n \tau}{(1-\eta) \hat{S}}\right)} = \frac{(1+\eta) \hat{\Omega}^2 \left(1-\frac{n\tau}{(1-\eta) \hat{S}}\right)+m \tau+\frac{(1+\eta) \hat{\Omega}^2 n \tau}{(1-\eta) \hat{S}}}{(1-\eta) \hat{S}\left(1- \frac{n \tau}{(1-\eta) \hat{S}}\right)}\nonumber\\
&=\left(\frac{1+\eta}{1-\eta}\right) \frac{\hat{\Omega}^2}{\hat{S}} + \frac{m \tau (1-\eta) \hat{S} + n \tau (1+\eta) \hat{\Omega}^2}{(1-\eta) \hat{S} \left((1-\eta) \hat{S}-n \tau \right)}.
\end{align}
We have 
\begin{align}\label{UpperBound2}
\frac{m \tau (1-\eta) \hat{S} + n \tau (1+\eta) \hat{\Omega}^2}{(1-\eta) \hat{S} \left((1-\eta) \hat{S}-n \tau \right)} \firstinequal \frac{m \tau (1-\eta) + n \tau (1+\eta)}{(1-\eta) \left((1-\eta) \hat{S}-n \tau \right)} \secondinequal \frac{m \tau (1-\eta) + n \tau (1+\eta)}{(1-\eta) \left((1-\eta) n s-n \tau \right)},    
\end{align}
where (a) follows from Lemma \ref{boundOmegaS2}, i.e., $\hat{\Omega}^2 \le \hat{S}$, and (b) follows since $S \ge n s$, for some $s > \frac{\tau}{1-\eta}$. 
For $n \gg m$, \eqref{UpperBound2} is approximated as 
\begin{align}\label{LowerBoundUpperApprox}
\frac{m \tau (1-\eta) \hat{S} + n \tau (1+\eta) \hat{\Omega}^2}{(1-\eta) \hat{S} \left((1-\eta) \hat{S}-n \tau \right)} \le \frac{m \tau (1-\eta) + n \tau (1+\eta)}{(1-\eta) \left((1-\eta) n s-n \tau \right)} \approx \frac{\tau}{(1-\eta)s - \tau}.
\end{align}
Plugging \eqref{LowerBoundUpperApprox} into \eqref{UpperBound1} yields 
\begin{align}
\frac{\hat{\Omega}^2}{\hat{S}} \le \left(\frac{1+\eta}{1-\eta}\right) \frac{\Omega^2}{S} + \frac{\tau}{(1-\eta)s - \tau}.    
\end{align}

\end{proof}

\section{Error bound on ratio of estimators}

\begin{lemma}\label{boundOmegaS2}
Assuming that $d_{\rm{max}}/d_{\rm{min}}^2 \le \frac{(n-1)}{2}$, where $d_{\rm{max}} = \max_{i,j} \left\| \boldsymbol{x}_i - \boldsymbol{x}_j \right\|^2$ and $d_{\rm{min}} = \min_{i,j, i \ne j} \left\| \boldsymbol{x}_i - \boldsymbol{x}_j \right\|^2$, we have $\Omega^2 \le S$.
\end{lemma}
\begin{proof}
We have
\begin{subequations}
\begin{align}
\Omega^2 &= \frac{1}{n} \sqrt{ 2 \operatorname{Tr}\left(\boldsymbol{Y}^{\top}\boldsymbol{L_ {X}{Y}}\right)}=\frac{1}{n} \sqrt{\sum_{i, j}{W^{X}_{i, j}} \left\| \boldsymbol{y}_i - \boldsymbol{y}_j \right\|^2} \nonumber \\
& = \frac{1}{n} \sqrt{ \sum_{i,j}{[\boldsymbol{JE_XJ}]_{i, j}} \left\| \boldsymbol{y}_i - \boldsymbol{y}_j \right\|^2},\\
S & = \frac{1}{n^4} \sum_{i,j} \left\| \boldsymbol{x}_i - \boldsymbol{x}_j \right\|^2 \sum_{i,j} \left\| \boldsymbol{y}_i - \boldsymbol{y}_j \right\|^2.
\end{align}
\end{subequations}
Assuming that there are a total of $c$ classes in the dataset, we denote by $\mathcal{C}_l$ the set of data samples indices that belong to the label $l$, $l =0, ..., c-1$.
For $i \in \mathcal{C}_l$, we assume that $\boldsymbol{y}_i$ is a vector with only one non-zero entry at the $l$-th coordinate, where for simplicity we set the value in the $l$-th coordinate to 1.
Accordingly, for $l=0, ..., c-1$, we have
\newcommand\firstequal{\mathrel{\overset{\makebox[0pt]{\mbox{\normalfont\tiny\sffamily (a)}}}{=}}}
\begin{align}
\left\| \boldsymbol{y}_i - \boldsymbol{y}_j \right\|^2 = 
\begin{cases}
0,& i, j \in \mathcal{C}_l \\
2,& \mbox{otherwise},
\end{cases}
\end{align}
which results in 
\begin{align}\label{OmegaIneq}
\Omega^2 &= \frac{1}{n} \sqrt{ 2\sum_{i,j, (i, j) \notin \mathcal{C}_l^2}{[\boldsymbol{JE_XJ}]_{i, j}}}    \nonumber\\
& \le \frac{1}{n} \sqrt{ 2\sum_{i,j}{[\boldsymbol{JE_XJ}]_{i, j}}} \nonumber\\
&\firstequal \frac{1}{n} \sqrt{\frac{2}{n} \sum_{i,j} \left\| \boldsymbol{x}_i - \boldsymbol{x}_j \right\|^2},
\end{align}
where (a) follows from Lemma \ref{LemmaJExJ}. 
For simplicity of the analysis, we assume that the number of data samples with each class is uniform, i.e., $\left| \mathcal{C}_l \right| = n/c$, in which case $\sum_{i,j} \left\| \boldsymbol{y}_i - \boldsymbol{y}_j \right\|^2 = 2 n^2 (1-1/c)$. 
Accordingly, we have
\begin{align}\label{S2Eq}
S &= \frac{2(1-1/c)}{n^2} \sum_{i,j} \left\| \boldsymbol{x}_i - \boldsymbol{x}_j \right\|^2, 
\end{align}
which is minimized for $c=2$ (assuming that $c \ge 2$), i.e.,
\begin{align}\label{S2Ineq}
S \ge \frac{1}{n^2} \sum_{i,j} \left\| \boldsymbol{x}_i - \boldsymbol{x}_j \right\|^2.
\end{align}
According to \eqref{OmegaIneq} and \eqref{S2Ineq}, in order to prove $\Omega^2 \le S$, it suffices to show that 
\begin{align}\label{SimplerIneq}
\sum_{i,j} \left\| \boldsymbol{x}_i - \boldsymbol{x}_j \right\|^2 \le \frac{1}{2n} \Big( \sum_{i,j} \left\| \boldsymbol{x}_i - \boldsymbol{x}_j \right\|^2 \Big)^2.    
\end{align}
We have
\begin{align}
\sum_{i,j} \left\| \boldsymbol{x}_i - \boldsymbol{x}_j \right\|^2 &\le n(n-1) d_{\rm{max}}, \nonumber\\
\frac{1}{2n} \Big( \sum_{i,j} \left\| \boldsymbol{x}_i - \boldsymbol{x}_j \right\|^2 \Big)^2 &\ge \frac{n^2(n-1)^2}{2n} d_{\rm{min}}^2.
\end{align}
Accordingly, having $d_{\rm{max}} \le \frac{(n-1)}{2}d_{\rm{min}}^2$ guarantees \eqref{SimplerIneq}.

\end{proof}






\nocite{langley00}



\section{Double-centering lemmas}
\begin{lemma}\label{LemmaJExJ}

For the centering matrix given by $\boldsymbol{J}=\boldsymbol{I}-\frac{1}{n} \boldsymbol{e e}^{\top}$, we have the following property when applied to Euclidean distance matrices of data sample $X$ denoted by $E_X$,  $$\sum_{i,j}\left[\boldsymbol{J E_{X} J}\right]_{i ,j} = \frac{1}{n} \sum_{i,j} d_{i,j}^2(X) $$
\end{lemma}

\begin{proof}
We have
\begin{align}
&\boldsymbol{J E_{X} J}=\left(\boldsymbol{I}-\frac{1}{n} \boldsymbol{e e}^{\top}\right) \boldsymbol{E_{X}}\left(\boldsymbol{I}-\frac{1}{n} \boldsymbol{e e}^{\top}\right)\\
&=\boldsymbol{E_{X}}-\frac{1}{n} \boldsymbol{e e}^{\top} \boldsymbol{E_{X}}-\frac{1}{n} \boldsymbol{E_{X} e e}^{\top}+\frac{1}{n^{2}} \boldsymbol{e e}^{\top} \boldsymbol{E_{X} e e}^{\top},
\end{align}
according to which
\begin{align}\label{JExJij}
&\left(\boldsymbol{J E_{X} J}\right)_{i j}=d_{i j}^{2}(\boldsymbol{X})-\frac{1}{n} \sum_{i'=1}^{n} d_{i' j}^{2}(\boldsymbol{X})-\frac{1}{n} \sum_{j'=1}^{n} d_{i j'}^{2}(\boldsymbol{X})+\frac{1}{n^{2}} \sum_{i' j'} d_{i' j'}^{2}(\boldsymbol{X}).
\end{align}
Summing \eqref{JExJij} over all $i,j$ yields
\begin{align}
\sum_{i, j}\left[\boldsymbol{J E_{X} J}\right]_{i j}=&\sum_{i, j} d_{i j}^{2}(\boldsymbol{X})-\frac{1}{n} \sum_{i, j} \sum_{i^{\prime}=1}^{n} d_{i^\prime j}^{2}(\boldsymbol{X})\\
&-\frac{1}{n} \sum_{i, j} \sum_{j^\prime=1}^{n} d_{i^{\prime} j^{\prime}}^{2}(\boldsymbol{X})+\frac{n(n-1)}{n^{2}} \sum_{i', j^{\prime}} d_{i^{\prime} j^{\prime}}^{2}(\boldsymbol{X})\\
=&\sum_{i j} d_{i j}^{2}(\boldsymbol{X})-\frac{n-1}{n} \sum_{i j} d_{i j}^{2}(\boldsymbol{X})-\frac{n-1}{n} \sum_{i, j} d_{ij}^{2}(\boldsymbol{X})+\frac{n-1}{n} \sum_{i, j} d_{i j}^{2}(\boldsymbol{X})\\
=&\frac{1}{n} \sum_{i, j} d_{i j}^{2}(\boldsymbol{X}).
\end{align}
\end{proof}

\begin{lemma}
The data matrix $\boldsymbol{X}$ and its Euclidean distance matrix $\boldsymbol{E_X}$ can be connected using the centering matrix given by $\boldsymbol{J}=\boldsymbol{I}-\frac{1}{n} \boldsymbol{e e}^{\top}$ as,
$$\boldsymbol{JX^TXJ} = -\frac{1}{2} \boldsymbol{JE_XJ}$$
\end{lemma}
\begin{proof}
\begin{align}
&\left(\left\|\bold{x}_{i}-\bold{x}_{j}\right\|^{2}=\bold{x}_{i}^{\top} \bold{x}_{i}+\bold{x}_{j}^{\top} \bold{x}_{j}-2 \bold{x}_{i}^{\top} \bold{x}_{j}\right)\\
&\bold{x}_{i}^{T} \bold{x}_{j}=-\frac{1}{2}\left(\left\|\bold{x}_{i}-\bold{x}_{j}\right\|^{2}-\left\|\bold{x}_{i}\right\|^{2}-\left\|\bold{x}_{j}\right\|^{2}\right)\\
&\bold{X}^T\bold{X} = \frac{-1}{2} \bold{E}_\bold{X} + \frac{1}{2}\bold{\delta.1}^T + \frac{1}{2}\bold{1 \delta}^T
\textrm{ where } \delta\textrm{ is a vector with }\delta_i=\|\bold{x}_i\|.
\end{align}
$$ \bold{JX^TXJ} = -\frac{1}{2} \bold{JE_XJ} + \frac{1}{2}J\bold{\delta}.\bold{1}^T\bold{J} + \frac{1}{2}\bold{J1 \delta}^T\bold{J}
$$
But since, $\bold{1}^T\bold{J} = 0$ and $\bold{J1}=\bold{0}$, we have $$\bold{JX}^T\bold{XJ} = -\frac{1}{2} \bold{JE_XJ}$$ 
\end{proof}

\section{Proof of Lemma \ref{boundOmegaS2}}
We have
\begin{align}
\Omega^2(\bold{X,Y}) &= \frac{1}{n} \sqrt{ 2 \operatorname{Tr}\left(\bold{Y}^{\top}\bold{L_ {X}{Y}}\right)}=\frac{2}{n^2} \sum_{i j}{\bold{W}^\bold{{X}}_{i j}} d_{i j}^{2}(\bold{Y}) \nonumber \\
& = \frac{1}{n} \sqrt{ \sum_{ij}{[\bold{JE_XJ}]_{i j}} d_{i j}^{2}(\bold{Y})}.\nonumber
\end{align}

\section{Laplacian formulation of\label{lapOmega}  $\Omega^2(\bold{X,Y})$ \cite{vepakomma2018supervised}} 
\begin{lemma}Distance covariance $\hat{\Omega}^2\boldsymbol{(X,Y)}$ can be estimated using the Euclidean distance matrix $\boldsymbol{E_X}$ formed over the rows in $\boldsymbol{X}$, the double-centering matrix $\boldsymbol{J=I}-n^{-1}\boldsymbol{ee^T}$, to form an adjacency matrix given by $\boldsymbol{W(X)=JE_XJ}$ and a corresponding graph Laplacian $\boldsymbol{L^W(X)=D(W(X))-W(X)}$ where $\boldsymbol{D(W(X))}$ is the degree matrix of $\boldsymbol{W(X)}$ to get $${\hat{\Omega}{^2}(\boldsymbol{X},\boldsymbol{Y}) = \frac{1}{n} \sqrt{ 2 \operatorname{Tr}\left(\boldsymbol{Y^{\top}L^W (X){Y}}\right)}\nonumber = \frac{1}{n} \sqrt{ 2 \operatorname{Tr}\left(\boldsymbol{X^{\top}L^W( Y){X}}\right)}} $$
\end{lemma}
\begin{proof}
This result and corresponding proof is from \cite{vepakomma2018supervised}. We replicate the same over here for quick reference. 
 Given matrices $\widehat{\boldsymbol{E}}_{\boldsymbol{X}}, \widehat{\boldsymbol{E}}_{\boldsymbol{Y}}$, and column centered matrices $\widetilde{\boldsymbol{X}}, \tilde{\boldsymbol{Y}}$, from result of \cite{torgerson1952multidimensional} we have that $\widehat{\boldsymbol{E}}_{\boldsymbol{X}}=-2 \widetilde{\boldsymbol{X}} \widetilde{\boldsymbol{X}}^{T}$ and $\widehat{\boldsymbol{E}}_{\boldsymbol{Y}}=-2 \widetilde{\boldsymbol{Y}} \tilde{\boldsymbol{Y}}^{T}$. In the problem of multidimensional scaling (MDS) (Borg and Groenen, 2005), we know for a given adjacency matrix say $\boldsymbol{W}$ and a Laplacian matrix $\boldsymbol{L}$,
$$
\operatorname{Tr}\left(\boldsymbol{X}^{T} \boldsymbol{L} \boldsymbol{X}\right)=\frac{1}{2} \sum_{i, j}[\boldsymbol{W}]_{i j}[\boldsymbol{E}_\boldsymbol{X}]_{i, j}
$$
Now for the Laplacian $\boldsymbol{L}=\boldsymbol{L}_{\boldsymbol{X}}$ and adjacency matrix $\boldsymbol{W}=\widehat{\boldsymbol{E}}_{\boldsymbol{Y}}$ we can represent $\operatorname{Tr}\left(\boldsymbol{X}^{T} \boldsymbol{L}_{\boldsymbol{Y}} \boldsymbol{X}\right)$ in terms of $\widehat{\boldsymbol{E}}_{\boldsymbol{Y}}$ as follows,
$$
\operatorname{Tr}\left(\boldsymbol{X}^{T} \boldsymbol{L}_{\boldsymbol{Y}} \boldsymbol{X}\right)=\frac{1}{2} \sum_{i, j=1}^{n}\left[\widehat{\boldsymbol{E}}_{\boldsymbol{Y}}\right]_{i, j}\left[\boldsymbol{E}_{\boldsymbol{X}}\right]_{i, j}
$$
From the fact $\left[\boldsymbol{E}_{\boldsymbol{X}}\right]_{i, j}=\left(\left\langle\widetilde{\boldsymbol{x}}_{i}, \widetilde{\boldsymbol{x}}_{i}\right\rangle+\left\langle\widetilde{\boldsymbol{x}}_{j}, \widetilde{\boldsymbol{x}}_{j}\right\rangle-2\left\langle\widetilde{\boldsymbol{x}}_{i}, \widetilde{\boldsymbol{x}}_{j}\right\rangle\right)$, and also $\widehat{\boldsymbol{E}}_{\boldsymbol{X}}=-2 \widetilde{\boldsymbol{X}} \widetilde{\boldsymbol{X}}^{T}$ we get
$$
\begin{aligned}
\operatorname{Tr}\left(\boldsymbol{X}^{T} \boldsymbol{L}_{\boldsymbol{Y}} \boldsymbol{X}\right)=&-\frac{1}{4} \sum_{i, j=1}^{n}\left[\widehat{\boldsymbol{E}}_{\boldsymbol{Y}}\right]_{i, j}\left(\left[\widehat{\boldsymbol{E}}_{\boldsymbol{X}}\right]_{i, i}+\left[\widehat{\boldsymbol{E}}_{\boldsymbol{X}}\right]_{j, j}-2\left[\widehat{\boldsymbol{E}}_{\boldsymbol{X}}\right]_{i, j}\right) \\
=& \frac{1}{2} \sum_{i, j}\left[\widehat{\boldsymbol{E}}_{\boldsymbol{X}}\right]_{i, j}\left[\widehat{\boldsymbol{E}}_{\boldsymbol{Y}}\right]_{i, j}-\frac{1}{4} \sum_{j}^{n}\left[\widehat{\boldsymbol{E}}_{\boldsymbol{X}}\right]_{j, j} \sum_{i}^{n}\left[\widehat{\boldsymbol{E}}_{\boldsymbol{Y}}\right]_{i, j} \\
&-\frac{1}{4} \sum_{i}^{n}\left[\widehat{\boldsymbol{E}}_{X}\right]_{i, i} \sum_{j}^{n}\left[\widehat{\boldsymbol{E}}_{\boldsymbol{Y}}\right]_{i, j}
\end{aligned}
$$
Since $\widehat{\boldsymbol{E}}_{\boldsymbol{X}}$ and $\widehat{\boldsymbol{E}}_{\boldsymbol{Y}}$ are double centered matrices $\sum_{i=1}^{n}\left[\widehat{\boldsymbol{E}}_{\boldsymbol{Y}}\right]_{i, j}=\sum_{j=1}^{n}\left[\widehat{\boldsymbol{E}}_{\boldsymbol{Y}}\right]_{i, j}=0$ it follows that
$$
\operatorname{Tr}\left(\boldsymbol{X}^{T} \boldsymbol{L}_{\boldsymbol{Y}} \boldsymbol{X}\right)=\frac{1}{2} \sum_{i, j}\left[\widehat{\boldsymbol{E}}_{\boldsymbol{X}}\right]_{i, j}\left[\widehat{\boldsymbol{E}}_{\boldsymbol{Y}}\right]_{i, j}
$$
It also follows that
$$
\hat{\nu^{2}}(\boldsymbol{X}, \boldsymbol{Y})=\frac{1}{n^{2}} \sum_{i, j=1}^{n}\left[\widehat{\boldsymbol{E}}_{\boldsymbol{Y}}\right]_{i, j}\left[\boldsymbol{E}_{\boldsymbol{X}}\right]_{i, j}=\frac{2}{n^{2}} \operatorname{Tr}\left(\boldsymbol{X}^{T} \boldsymbol{L}_{\boldsymbol{Y}} \boldsymbol{X}\right)
$$
Similarly, we can express the sample distance covariance using Laplacians $\boldsymbol{L}_{\boldsymbol{X}}$ and $\boldsymbol{L}_{\boldsymbol{Y}}$ as
$$
\hat{\nu^{2}}(\boldsymbol{X}, \boldsymbol{Y})=\left(\frac{2}{n^{2}}\right) \operatorname{Tr}\left(\boldsymbol{X}^{T} \boldsymbol{L}_{\boldsymbol{Y}} \boldsymbol{X}\right)=\left(\frac{2}{n^{2}}\right) \operatorname{Tr}\left(\boldsymbol{Y}^{T} \boldsymbol{L}_{\boldsymbol{X}} \boldsymbol{Y}\right)
$$
The sample distance variances can be expressed as $\hat{\nu}^{2}(\boldsymbol{X}, \boldsymbol{X})=\left(\frac{2}{n^{2}}\right) \operatorname{Tr}\left(\boldsymbol{X}^{T} \boldsymbol{L}_{\boldsymbol{X}} \boldsymbol{X}\right)$ and $\hat{\nu}^{2}(\boldsymbol{Y}, \boldsymbol{Y})=\left(\frac{2}{n^{2}}\right) \operatorname{Tr}\left(\boldsymbol{Y}^{T} \boldsymbol{L}_{\boldsymbol{Y}} \boldsymbol{Y}\right)$ substituting back into expression of sample
\end{proof}

\section{$S(\bold{X,Y})$ as a sum of directional variances \label{sdv}}
\begin{lemma}\vspace{-2mm}We now express the denominator term of the test-statistic as a specific sum of directional variances as follows.
The denominator in the test-statistic, $${\hat{S}} (\boldsymbol{X,Y})=\frac{1}{n^{2}} \sum_{k, l=1}^{n} \left\| \boldsymbol{x}_{k}-\boldsymbol{x}_{l}\right\|^{2} \frac{1}{n^{2}} \sum_{k, l=1}^{n}\left\| \boldsymbol{y}_{k}-\boldsymbol{y}_{l}\right\|^{2},$$ can be expressed as a sum of directional variances as $\hat{S}(\boldsymbol{X,Y})=\sum_i\phi_{\boldsymbol{X}}(\boldsymbol{g_i})$ where $\boldsymbol{g}_i$ is the $i$'th column of a matrix $\boldsymbol{G}$ such that  $\boldsymbol{L^S}=\boldsymbol{G}\boldsymbol{G}^T$ for a graph Laplacian matrix $\bold{L^{S}}$ given by $\bold{L^{S}}=n\boldsymbol{I-ee^T}$. We use the superscript $\boldsymbol{S}$ to distinguish from the Laplacian $\boldsymbol{L^W}$ used in the expression for $\hat{\Omega}^2$.\end{lemma}
\vspace{-2mm}
\begin{proof}
We first start simple, by expressing the Euclidean distance between any two pairs of points using a matrix trace formulation as follows.
$$
d_{ij}^{2}(\boldsymbol{X}) =\sum_{a=1}^{m} \boldsymbol{x}_{a}^{\top}\left(\boldsymbol{e}_{i}-\boldsymbol{e}_{j}\right)\left(\boldsymbol{e}_{i}-\boldsymbol{e}_{j}\right)^{\top} \boldsymbol{x}_{a} 
=\sum_{a=1}^{m} \boldsymbol{x}_{a}^{\top} \boldsymbol{A}_{ij} \boldsymbol{x}_{a} 
=\operatorname{tr} \boldsymbol{X}^{\top} \boldsymbol{A}_{ij} \boldsymbol{X} 
$$
where $\boldsymbol{A}_{ij} = (e_i-e_j)(e_i-e_j)^{\top}$. Similarly, $
\sum_{i < j} d_{i j}(\boldsymbol{X})=\operatorname{tr} \boldsymbol{X}^{\top}\left(\sum_{i<j}d_{i j}^{-1}(\boldsymbol{X}) \boldsymbol{A}_{i j}\right) \boldsymbol{X}
$.
We now extend this to express the sum of all pairs of Euclidean distance matrices as follows
$$\eta(\boldsymbol{X}) =\sum_{i<j} w_{i j} d_{i j}^{2}(\boldsymbol{X})=\operatorname{tr} \boldsymbol{X}^{\top}\left(\sum_{i<j} w_{i j} \boldsymbol{A}_{i j}\right) \boldsymbol{X} 
=\operatorname{tr} \boldsymbol{X}^{\top} \boldsymbol{L^{S}} \boldsymbol{X}$$ where $\bold{L^S}=\left(\sum_{i<j} w_{i j} \boldsymbol{A}_{i j}\right)$. 
As we would like to express the sum of all pairs of distances, we consider the case where $\boldsymbol{W}_{ij}=1, \forall i,j \in [n]$.  Note that this matrix exactly has the structure of a graph Laplacian. This results in a graph Laplacian $\boldsymbol{L^S}$ where

$$L^{S}=\left[\begin{array}{ccc}
n-1 & -1 & -1 \\
-1 & \ddots& \vdots \\
-1 & \cdots & n-1
\end{array}\right]=n\boldsymbol{I-ee^T}$$
Therefore,
$$\boldsymbol{S^\alpha(X,Y)}=\eta(\boldsymbol{X})\eta(\boldsymbol{Y})=\frac{4  }{n^{4}}Tr\left(\boldsymbol{X}^{\top} \boldsymbol{L^{S} X}\right) \cdot \operatorname{Tr}\left(\boldsymbol{Y}^{T} \boldsymbol{L^{S} Y}\right)$$
Note that we could also express $Tr\left(\boldsymbol{X}^{\top} \boldsymbol{L^{S} X}\right)$ as  $$Tr\left(\boldsymbol{X}^{\top} \boldsymbol{L^{S} X}\right) =Tr\left(\boldsymbol{X}^{\top} \boldsymbol{GG^{T} X}\right) =Tr\left(\boldsymbol{G^{T} XX}^{\top} \boldsymbol{G}\right) =\sum_{i}^{k\leq n}\left(\boldsymbol{g_i}^{\boldsymbol{T}} \boldsymbol{X} \boldsymbol{X}^{\boldsymbol{T}}\boldsymbol{g_i}\right)=\sum_i\phi_{\boldsymbol{X}}(\boldsymbol{g_i}),$$ where $k$ is the rank of the matrix $\boldsymbol{L^{S}}$.
\end{proof}
\section{Population distance correlation}\label{pdc}
Using the definition of distance covariance, we have the following expression for the square of distance correlation \cite{szekely2007measuring} between random variables \[
	\left\{ \begin{array}{cc}
				  \frac{{\Omega^2}(\boldsymbol{X},\boldsymbol{Y})}{\sqrt{{\Omega^2}(\boldsymbol{X},\boldsymbol{X}){\Omega^2}(\boldsymbol{Y},\boldsymbol{Y})}}, 
				 & \text{if }{\Omega^2}(\boldsymbol{X},\boldsymbol{X}){\Omega^2}(\boldsymbol{Y},\boldsymbol{Y})>0\\
	        0, & \text{if } {\Omega^2}(\boldsymbol{X},\boldsymbol{X}){\Omega^2}(\boldsymbol{Y},\boldsymbol{Y})=0
	        \end{array} \right.
	\]

This always lies within the interval $[0,1]$ with $0$ referring to independence and $1$ referring to dependence.

\section{Families of some dependency measures}
\label{appdx:Families}
\begin{enumerate}
    \item \textbf{Energy statistics:} Distance correlation (DCOR), Brownian distance covariance, Partial distance correlation, Partial martingale difference correlation
    \item \textbf{Kernel covariance operators:} Constrained covariance (COCO), Hilbert-Schmidt independence criterion (HSIC), Kernel Target Alignment (KTA)
    \item \textbf{Integral probability  metrics:} Maximum mean discrepancey (MMD), Wasserstein distance,  Dudley metric and Fortet-Mourier metric, Total variation distance.
    \item \textbf{Information theoretic measures:} Mutual information, f-divergence, Renyi divergence, Hellinger distance, Total variation distance, Maximal information coefficient, Total information coefficient
\end{enumerate}
\bibliography{example_paper}